\documentclass{cedram-aif}

% ...

\usepackage{amsmath,amsfonts,amssymb}
\usepackage{graphicx}

\newtheorem{ethm}{Theorem}[section]

\newtheorem{ecor}[ethm]{Corollary}

\newtheorem{eprop}[ethm]{Proposition}

\newtheorem{elem}[ethm]{Lemma}

\newcommand{\proofend}{~$\rhd$}
\newcommand{\proofbegin}{~$\lhd$}

\newenvironment{eproof}
               {\noindent {\emph{\textbf{Proof}}}\\\proofbegin~}
               {\proofend\\}

\newtheorem{erem}[ethm]{Remark}

\newcommand{\dR}{\ensuremath{\mathbb{R}}}

\newcommand{\R}{\dR}
\newcommand{\p}[4]{{#3}\!\left#1{#4}\right#2}
\newcommand{\PAR}[1]{\ensuremath{{\left(#1\right)}}} % (1)
\newcommand{\BRA}[1]{\ensuremath{{\left\{#1\right\}}}} % {1}

\def\disp{\displaystyle}
\newcommand{\entf}[1]{{\rm{Ent}}_{#1}}
\newcommand{\ent}[2]{\p(){\entf{#1}}{#2}}
\newcommand{\bigO}[1]{\ensuremath{\mathop{}\mathopen{}O\mathopen{}\left(#1\right)}}
\newcommand{\smallO}[1]{\ensuremath{\mathop{}\mathopen{}o\mathopen{}\left(#1\right)}}

% ...

\title
[The Li-Yau inequality and applications]
{The Li-Yau inequality and applications  under a curvature-dimension condition}

\alttitle{In\'egalit\'e de Li-Yau et applications sous une condition de courbure-dimension}

\author{\firstname{Dominique} \lastname{Bakry}}

\address{Institut de math\'ematiques de Toulouse,\\
 Umr Cnrs 5219,\\
  Universit\'e Toulouse 3. }

\email{bakry@math.univ-toulouse.fr}
\thanks{This research was supported  by the French ANR-12-BS01-0019 STAB project.}

\author{\firstname{Fran{\c c}ois} \lastname{Bolley}}
\address{Laboratoire de probabilit\'es et mod\`eles al\'eatoires,\\
 Umr Cnrs 7599,\\
 Universit\'e Paris 6. francois.bolley@upmc.fr}
\email{francois.bolley@upmc.fr}

\author{\firstname{Ivan} \lastname{Gentil}}
\address{Institut Camille Jordan,\\
 Umr Cnrs 5208,\\
  Universit\'e Claude Bernard Lyon 1}
\email{gentil@math.univ-lyon1.fr}
\keywords{Li-Yau inequality, Harnack inequality, heat kernel bounds, Ricci curvature. }
  
\altkeywords{in\'egalit\'e de Li-Yau, in\'egalit\'e de Harnack, noyaux de la chaleurs, courbure de Ricci}

\subjclass{58J35,46-XX, 60H15}

\begin{document}
%% Abstract 
\begin{abstract}
  We prove a global  Li-Yau inequality for a general Markov semigroup under a curvature-dimension condition. This inequality is stronger than all classical Li-Yau type inequalities known to us.  On a Riemannian manifold, it is equivalent  to a new parabolic Harnack inequality, both in negative and positive curvature, giving new subsequent bounds on the heat kernel of the semigroup. Under positive curvature we moreover reach ultracontractive bounds by a direct and robust method.  
\end{abstract}

%% French abstract
\begin{altabstract}
Nous obtenons une in\'egalit\'e de type Li-Yau pour un semi-groupe de Markov g\'en\'eral, sous une condition de courbure-dimension. A notre connaissance, cette nouvelle in\'egalit\'e renforce toutes les in\'egalit\'es de ce type. Sur une vari\'et\'e riemannienne, elle est \'equivalente \`a une nouvelle in\'egalit\'e de Harnack parabolique, en courbure positive ou n\'egative, et induit des bornes pertinentes sur le noyau de la chaleur associ\'e. En courbure positive, elle permet d'atteindre des bornes ultracontractives par une m\'ethode directe et robuste.
\end{altabstract}

\maketitle

\section{Introduction}

In their seminal paper~\cite{li-yau}, P. Li and S.-T. Yau proved that on a Riemannian manifold $M$ with dimension $n$ and non-negative Ricci curvature, for any positive function $f$ and any $t>0$,
\begin{equation}
\label{eq-li-yau}
-\Delta(\log P_t f)\leq \frac{n}{2t},
\end{equation}
where $\Delta$ is the Laplace-Beltrami operator on $M$ and $(t,x)\mapsto P_tf(x)$ is the solution to the heat equation $\partial_t u=\Delta u$ with initial condition $f$. Equivalently, the Li-Yau inequality can be written~: 
$$
\frac{|\nabla P_tf|^2}{(P_t f)^2}\leq \frac{\Delta P_t f}{P_t f}+\frac{n}{2t}, 
$$
where $|\nabla P_tf|$ stands for the length of $\nabla P_tf$. 
This gradient estimate is a crucial step towards parabolic Harnack inequalities and various subsequent on and off-diagonal bounds on heat kernels. It is optimal since the equality is achieved for the heat kernel on the Euclidean space, that is, when $f$ converges to a Dirac mass.

This inequality has been generalised to Riemannian manifolds with a Ricci curvature bounded from below by a real constant $\rho$. Still in~\cite{li-yau}, P.~Li and S.-T. Yau proved that if the Ricci curvature is bounded from below by $\rho=-K < 0$ then for any $\alpha>1$
$$
\frac{|\nabla P_tf|^2}{(P_t f)^2}\leq \alpha \frac{\Delta P_t f}{P_t f}+\frac{n\alpha^2K}{2(\alpha-1)}+\frac{n\alpha^2}{2t}. 
$$
Improving this inequality  towards a form which would be optimal for small, large as well as intermediate times has been the goal of many subsequent works. In this direction and still in the negative curvature case, let us mention the various bounds
\begin{equation}\label{eq:davies}
\frac{|\nabla P_tf|^2}{(P_t f)^2}\leq \alpha \frac{\Delta P_t f}{P_t f}+\frac{n\alpha^2K}{4(\alpha-1)}+\frac{n\alpha^2}{2t} 
\end{equation}
for any $\alpha>1$, derived by B. Davies in~\cite[Chapter 5.3]{davies},  and 
$$
\frac{|\nabla P_tf|^2}{(P_t f)^2} - \frac{\Delta P_t f}{P_t f} \leq \sqrt{2nK} \, \sqrt{ \frac{\vert \nabla P_t f \vert^2}{P_t f}+\frac{n}{2t}+2nK} + \frac{n}{2t}
$$
derived by S.-T. Yau~\cite{yau1}, itself improved by the first author and Z. Qian~\cite{bakry-qian}, as
\begin{equation}\label{eq:bakryqian}
\frac{|\nabla P_tf|^2}{(P_t f)^2} - \frac{\Delta P_t f}{P_t f} \leq \sqrt{nK} \, \sqrt{ \frac{\vert \nabla P_t f \vert^2}{P_t f}+\frac{n}{2t}+\frac{nK}{4}} + \frac{n}{2t}. 
\end{equation}
Meanwhile, R. Hamilton~\cite{hamilton} had proved the inequality
\begin{equation}\label{eq:hamilton}
\frac{|\nabla P_tf|^2}{(P_t f)^2} \leq e^{2Kt} \frac{\Delta P_t f}{P_t} + \frac{n}{2t} e^{4Kt} 
\end{equation} 
and, most recently,  J.~Li and X.~Xu~\cite{li-xu} have obtained the bound
\begin{equation}\label{eq:lixu}
\frac{\vert \nabla P_tf \vert^2}{(P_tf)^2}\leq \left(1+ \frac{\sinh(2 Kt)-2Kt}{2\sinh^2(Kt)}\right) \,  \frac{\Delta P_tf}{P_tf}+\frac{nK}{2} \big( 1+\coth(Kt) \big).
 \end{equation}
All these inequalities are based on the maximum principle and  inequalities \eqref{eq:davies} to \eqref{eq:lixu} are not comparable to each other and none of them reaches an optimal form. 

There is a rich literature on extensions of the Li-Yau inequality~\eqref{eq-li-yau} to diverse settings and evolution equations, currently almost 400 citations on MathSciNet or Zentralblatt. Let us mention in particular the most recent~\cite{souplet-zhang,wang2010,baudoin-garofalo,qian2,lee2013,qian-zhang-zhu,garofalo-mondino,qian}.   

\bigskip

In this article we shall prove the following general Li-Yau inequality,
{\it 
in negative and positive curvature. To our knowledge it improves  all Li-Yau inequalities obtained so far.  }  Assume that the Ricci curvature of the $n$-dimensional Riemannian manifold is (uniformly) bounded from below by $\rho\in\R^*$. Then, see Corollary~\ref{cor-li-yau-gene},
$$
\frac{4}{n\rho}\frac{\Delta P_tf}{P_tf} <  1+\frac{\pi^2}{\rho^2t^2}
$$ 
and 
$$
\frac{|\nabla P_t f|^2}{(P_tf)^2} < \frac{n}{2}\Phi_t\left(\frac{4}{n\rho}\frac{\Delta P_tf}{P_tf}\right)
$$
for any positive function $f$ and $t>0$, where 
$$
\Phi_t(x)=\left\{
\begin{array}{ll}
\disp \frac{\rho}{2}\PAR{x-2+2\sqrt{1-x}\coth(\rho t \sqrt{1-x})}, & x\leq 1\\
\disp \frac{\rho}{2}\PAR{x-2+2\sqrt{x-1}\cot(\rho t \sqrt{x-1})}, & 1\leq x<1+ \displaystyle \frac{\pi^2}{\rho^2t^2}.
\end{array}
\right.
$$

This result will be obtained as an extension to any curvature lower bound of the equivalence between the following properties, which is due to~\cite{bakry-ledoux}:  
\begin{enumerate}
\item The Ricci curvature of $M$ is non-negative. 
\item For any smooth positive function $f$ and $t>0$, 
\begin{equation}
\label{li-yau-bl}
\exp\left(-\frac{2}{nP_t f}(\ent{P_t}{f}+t\Delta P_tf \right)\leq 1+\frac{2t}{nP_t f}\PAR{\Delta P_tf-\frac{|\nabla P_tf|^2}{P_tf}}
\end{equation}
where $\ent{P_t}{f}=P_t(f\log f)-P_t f \; \log P_tf$. 
\item For any smooth positive function $f$ and $t>0$,
$$
\exp\left(\frac{2}{nP_t f}(\ent{P_t}{f}-t\Delta P_tf \right)\leq 1+\frac{2t}{nP_t f}\left(P_t\left(\frac{|\nabla f|^2}{f}\right)-\Delta P_tf\right).
$$
\item For any smooth positive function $f$ and $t>0$,
\begin{equation}
\label{commutation-positive}
P_tf \; \Delta (\log P_tf)\geq P_t(f\Delta \log f)\left(1+\frac{2t}{n}\Delta (\log P_t f)\right).
\end{equation}
\end{enumerate}

Inequality~\eqref{li-yau-bl} clearly implies the  Li-Yau inequality~\eqref{eq-li-yau}. Inequality~\eqref{commutation-positive} can be reformulated as a bound on the gradient of the heat kernel (commutation inequality), taking the dimension into account.  

\bigskip

The present paper gives a generalisation of this equivalence with any lower bound on the Ricci curvature,  positive or negative, instead only of a  non-negative one. As in~\cite{bakry-ledoux}, the result will be stated for general Markov semigroups, including diffusion semigroups on weighted Riemannian manifolds. 

\bigskip

 This work is organised as follows. In the next section we state this generalisation, for a Markov 
 diffusion semigroup under a  $CD(\rho,n)$ curvature-dimension condition. We also derive first consequences, including our main result: {\it a new Li-Yau inequality} under this curvature-dimension condition. The proof, which is very short and simple, is given in Section~\ref{sec-preuve}. Section~\ref{sec-applications} is devoted to applications: {\it ultracontractive bounds} in the positive curvature case, and {\it new Harnack inequalities} in the positive and negative cases, which are equivalent to our Li-Yau equality.  Finally,  in Section~\ref{sec-comparaison} we prove that our inequality implies all classical  Li-Yau inequalities known to us.

\section{Main result}
\label{sec-result}

\subsection{Markov triple and curvature-dimension condition}
\label{sec-gamma2}
A Markov diffusion triple $(E,\mu,\Gamma)$, as defined  in~\cite[Chapter 3]{bgl-book},  consists in a nice state space $E$ equipped with a Markov diffusion semigroup $(P_t)_{t\geq0}$ with infinitesimal generator $L$, carr\'e du champ $\Gamma$ and invariant and reversible $\sigma$-finite measure $\mu$.  The carr\'e du champ and $\Gamma_2$ operators are pointwise defined from the generator $L$ by 
$$
\Gamma(f,g)=\frac{1}{2}(L(fg)-fLg-gLf) \qquad \textrm{and} \qquad \Gamma_2 (f) = \frac{1}{2} (L \Gamma(f,f) - 2 \Gamma(f, Lf)) 
$$
for functions $f$ and $g$ in a suitable algebra $\mathcal A$ of functions from $E$ to $\R$. We let $\Gamma(f)=\Gamma(f,f)$. The generator $L$ is assumed to satisfy the diffusion property, that is, for any smooth function $\phi$ and $f\in\mathcal A$, 
$$
L\phi(f)=\phi'(f)Lf+\phi''(f)\Gamma(f).
$$

In the Markov triple setting, the abstract curvature-dimension condition $CD(\rho,n)$, for $\rho\in\R$ and $n\geq 1$, is satisfied when
$$
\Gamma_2 (f) \geq \rho \, \Gamma(f)+\frac{1}{n}(Lf)^2
$$
for any $f\in\mathcal A$.

Let us recall that under a $CD(\rho, n)$ condition with $\rho>0$, then the semigroup is ergodic, that is, $P_tf\rightarrow\int fd\mu$ in $L^2(\mu)$.

\bigskip

The main example of a Markov diffusion triple is a smooth, connected and complete weighted Riemannian manifold $(M,g)$ equipped with the generator $L=\Delta+\nabla V\cdot\nabla$, where  $\Delta$ is the Laplace-Beltrami operator and $V$ a smooth function on $M$, and with the measure $d\mu=e^Vd x$ where $d x$ is the Riemannian measure. In this case, the carr\'e du champ operator is  $\Gamma(f)=|\nabla f|^2$ where $ |\nabla f|$ is the length of the vector $\nabla f$, and the algebra of functions $\mathcal A$ consists in smooth and bounded functions on $M$. 

For instance, when $V=0$, the Bochner-Lichnerowicz inequality implies that the condition $CD(\rho,n)$ is satisfied when $M$ is a $n$-dimensional Riemannian manifold with a  Ricci curvature $Ric$ (uniformly)  bounded from below by $\rho$. For a general $V$, the condition holds on the $m$-dimensional manifold $M$ as soon as $m<n$ and
$$
Ric - \nabla^2 V \geq \rho + \frac{\nabla V \otimes \nabla V}{n-m} \cdot
$$

\bigskip

For a positive function $f$ on $E$ we let $\ent{P_t}{f} = P_t(f \log f) - P_t f \, \log P_t f$ and $\ent{\mu}{f} = \int f \log f \, d\mu - \int fd\mu \; \log \int fd\mu$.

\begin{erem}
In this work we shall deal with a symmetric Markov semigroup for convenience, but all the results proved here can be stated in a non-symmetric case. 
\end{erem}

\subsection{Li-Yau inequality under the $CD(\rho,n)$ condition}

\begin{ethm}[Local logarithmic Sobolev inequalities]
\label{theo-main}
Let $(E,\mu,\Gamma)$ be a Markov diffusion triple, $\rho\in\R^*$ and $n \geq 1$. Given a positive function $f$ on $E$ and $
t>0$, we set 
$$
X=\frac{4}{n\rho}\frac{LP_tf}{P_tf}
$$
and given $t>0$ we define the functions $\Phi_t$ and  $\tilde{\Phi}_t$ by  
\begin{equation}
\label{eq-defgh}
\Phi_t(x) =
\left\{ \begin{array}{ll}
\disp \frac{\rho}{2}\PAR{x-2+2\sqrt{1-x}\coth(\rho t \sqrt{1-x})}, & x\leq 1\\
\disp \frac{\rho}{2}\PAR{x-2+2\sqrt{x-1}\cot(\rho t \sqrt{x-1})}, & 1\leq x<1+ \displaystyle \frac{\pi^2}{\rho^2t^2}
\end{array}
\right.
\end{equation}
and $\tilde{\Phi}_t (x) = \Phi_t (x) -\rho x +2\rho $.
 Then the following properties are equivalent:
\begin{enumerate}
\item The Markov triple satisfies the $CD(\rho,n)$ condition. 
\item For any positive function $f\in\mathcal A$ and $t>0$,  then  $X< 1+ \displaystyle \frac{\pi^2}{\rho^2 t^2}$
and
\begin{multline}
\label{eq-logsob-inverse}
\exp\left(-\frac{2}{n}\frac{\ent{P_t}{f}}{P_t f}+\frac{t\rho}{2}X-\rho t\right)\leq\\
\left\{
\begin{array}{ll}
\disp\frac{\sinh(\rho t\sqrt{1-X})}{\rho\sqrt{1-X}}\PAR{ \Phi_t(X)-\frac{2}{n}\frac{\Gamma(P_tf)}{(P_tf)^2}}\,\, & {\rm{if}}\,\, X\leq 1\\
\disp\frac{\sin(\rho t\sqrt{X-1})}{\rho\sqrt{X-1}}\PAR{ \Phi_t(X)-\frac{2}{n}\frac{\Gamma(P_tf)}{(P_tf)^2}}\,\, & {\rm{if}}\,\, 1\leq X< 1+ \displaystyle \frac{\pi^2}{\rho^2 t^2}.\\
\end{array}
\right.
\end{multline}
\item For any positive function $f\in\mathcal A$ and $t>0$,  then $X< 1+ \displaystyle \frac{\pi^2}{\rho^2 t^2}$
and 
\begin{multline}
\label{eq-logsob}
\exp\left(\frac{2}{n}\frac{\ent{P_t}{f}}{P_t f}-\frac{t\rho}{2}X+\rho t\right)\leq\\
\left\{
\begin{array}{ll}
\disp\frac{\sinh(\rho t\sqrt{1-X})}{\rho\sqrt{1-X}}\PAR{\tilde{\Phi}_t(X) + \frac{2}{n} \frac{P_t(\Gamma(f)/f)}{P_t f}}\,\,&  {\rm{if}}\,\, X\leq 1\\
\disp\frac{\sin(\rho t\sqrt{X-1})}{\rho\sqrt{X-1}}\PAR{\tilde{\Phi}_t(X) + \frac{2}{n} \frac{P_t(\Gamma(f)/f)}{P_t f}}\,\, & {\rm{if}}\,\, 1\leq X< 1+ \displaystyle \frac{\pi^2}{\rho^2 t^2}.\\
\end{array}
\right.
\end{multline}

\end{enumerate}
\end{ethm}

As is \cite{bgl-book}, inequality  \eqref{eq-logsob} may be  called a local logarithmic Sobolev inequality and \eqref{eq-logsob-inverse} a local reverse logarithmic Sobolev inequality.
We observe that the right-hand sides of \eqref{eq-defgh}, \eqref{eq-logsob-inverse} and \eqref{eq-logsob} are continuous in $X=1$, justifying the above formulations.

\begin{ecor}
\label{cor-esti2}
Under the $CD(\rho,n)$ condition, for any positive function $f\in\mathcal A$ and $t>0$
\begin{equation}
\label{eq-esti2}
\frac{4}{n\rho}\frac{LP_t f}{P_tf}< 1+ \frac{\pi^2}{\rho^2 t^2}. 
\end{equation}
\end{ecor}

The term in  the right-hand side of~\eqref{eq-logsob-inverse} has to be positive, giving the following improved Li-Yau inequality:
\begin{ecor}[Improved Li-Yau inequality]
\label{cor-li-yau-gene}
For any Markov diffusion triple satisfying a $CD(\rho,n)$ condition with $\rho\in\R^*$ and $n\geq 1$, then 
\begin{equation}
\label{eq-li-yau-gene}
\frac{\Gamma( P_t f)}{(P_tf)^2}< \frac{n}{2} \, \Phi_t\left(\frac{4}{n\rho}\frac{L P_tf}{P_tf}\right)
\end{equation}
for any positive function $f\in\mathcal A$ and $t>0$, where the function $\Phi_t$ is defined in~\eqref{eq-defgh}. 
\end{ecor}

\begin{erem}
\begin{enumerate}
\item When $\rho$ tends to $0$ in Theorem~\ref{theo-main}, we exactly recover the estimates ii. and iii. given in the introduction in the case $\rho =0$. In particular the general Li-Yau inequality~\eqref{eq-li-yau-gene} converges  to the classical Li-Yau inequality~\eqref{eq-li-yau}. 
\item Combining the two inequalities~\eqref{eq-logsob-inverse} and~\eqref{eq-logsob} leads to a commutation type inequality similar to~\eqref{commutation-positive}, and which converges to~\eqref{commutation-positive} when $\rho$ tends to $0$. The inequality obtained in this general case where $\rho \neq 0$ is still equivalent to the $CD(\rho, n)$ condition, but is less appealing than~\eqref{commutation-positive}, which is why we omit it.  
\end{enumerate}
\end{erem}

\begin{erem}
We shall see in Section \ref{sec-comparaison} that for any $t$ the bound~\eqref{eq-li-yau-gene} improves upon all classical bounds recalled in the introduction.
\end{erem}

\begin{erem}
Our Li-Yau inequality~\eqref{eq-li-yau-gene} also improves on the following results obtained in \cite[Thms. 1 and 2]{bakry-qian}:

$\bullet$ Assume that $\rho >0$. Then \eqref{eq-li-yau-gene} holds for $X \leq X_0 :=1 +\pi^2/64$, where as above $X = 4 L P_t f / (n \rho P_t f)$; if $X > X_0$, then the bound holds with $\Phi_t(X)$ replaced by the tangent $\Phi_t(X_0) + \Phi_t'(X_0) (X-X_0)$. 

$\bullet$ Assume that $\rho<0$. Then \eqref{eq-li-yau-gene} holds for $X \leq 1$; moreover, for any $X$, then the bound holds with $\Phi_t(X)$ replaced by the tangent $\Phi_t(X_1) + \Phi_t'(X_1) (X-X_1)$ for any $X_1 \leq 1$.

But, by Lemma~\ref{lem-prop-phi} below and for any given $t>0$, the function $\Phi_t$ is a $\mathcal C^\infty$ and strictly concave function on the whole interval $(-\infty,1+\frac{\pi^2}{\rho^2t^2})$;  hence its graph is below its tangents, and \eqref{eq-li-yau-gene} improves upon these bounds in~\cite{bakry-qian}. 
\end{erem}

\begin{elem}[Properties of $\Phi_t$]
\label{lem-prop-phi}
For any $t>0$ and $\rho\in\R^*$, the function $\Phi_t$ is $\mathcal C^\infty$ and strictly concave on the interval $(-\infty, 1+\frac{\pi^2}{\rho^2 t^2})$. 
\end{elem}
\begin{eproof}
We have already observed that $\Phi_t$ is continuous in $X=1$. Moreover $\Phi_t$ is $\mathcal C^\infty$ since in $X=1$, the Taylor expansions (of any orders) are the same for $X=1^-$ or $X=1^+$.

Let us now prove that $\Phi_t$ is strictly concave, for instance in the case $\rho>0$, the case $\rho<0$ being similar.  For fixed $u\in(0,1)$, by direct computation, the  map $y \mapsto \frac{\sin(yu)}{\sin(y)}$ is increasing and positive on  $(0,\pi)$. This implies that the map 
$$
I: y \mapsto \int_0^1\PAR{\frac{\sin(yu)}{\sin(y)}}^2du=\frac{2y-\sin(2y)}{4y\sin^2(y)}
$$
is increasing on $(0,\pi)$. Hence the function $\Phi_t$ has a decreasing derivative
$$
\Phi_t'(x)=\rho/2-\rho^2 t \,  I(\rho t\sqrt{x-1})
$$
on  $(1,1+\frac{\pi^2}{\rho^2 t^2})$, so is strictly concave on this interval. 
The same argument can be performed on~$X<1$, starting from the function $y\mapsto \frac{\sinh(yu)}{\sinh(y)}$. This concludes the proof by differentiability of $\Phi_t$ in $X=1$.
\end{eproof}

\begin{figure}[!h]
\begin{center}
\includegraphics[width=9cm]{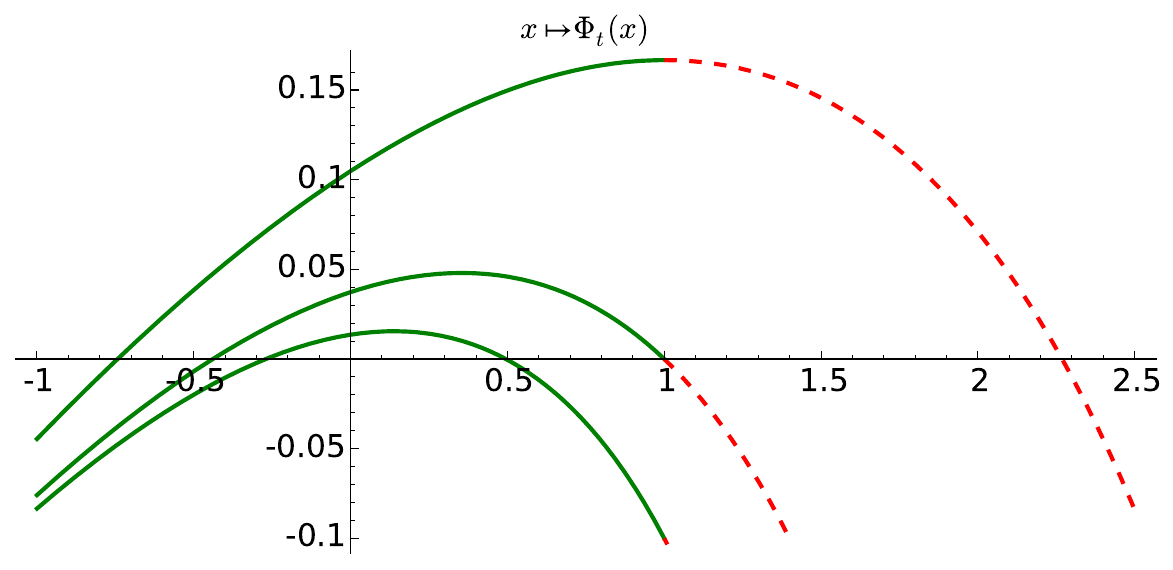}
\caption{Graphs of $\Phi_t$  for positive curvatures}
\label{fig-1}
\end{center}
\end{figure}

\begin{figure}[!h]
\begin{center}
\includegraphics[width=9cm]{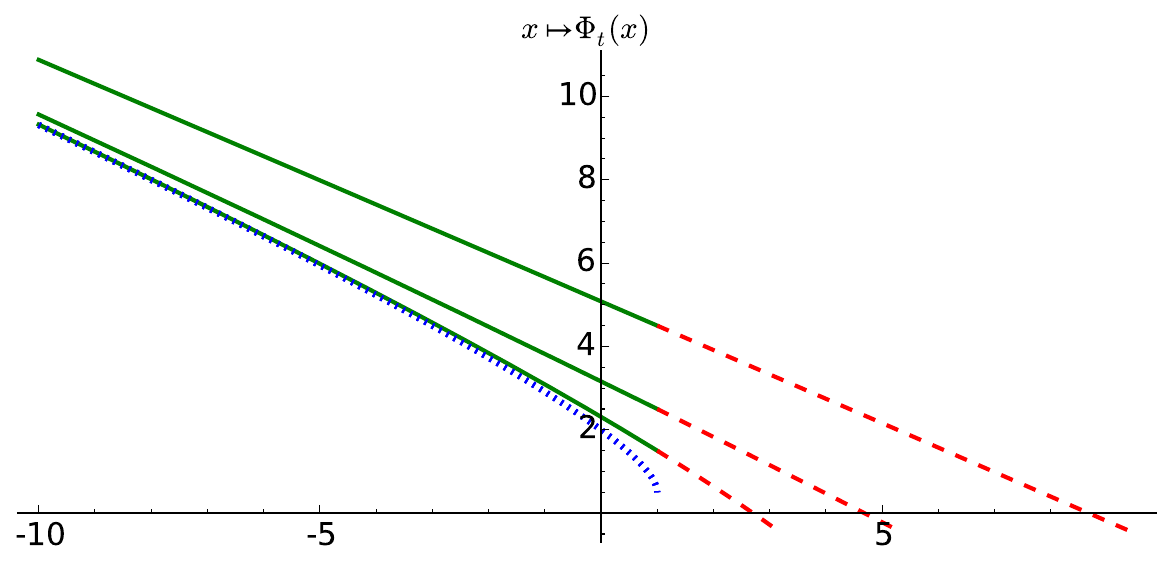}
\caption{Graphs of $\Phi_t$ for  negative curvatures}
\label{fig-2}
\end{center}
\end{figure}

 On the figure~\ref{fig-1} (resp. on the figure~\ref{fig-2}), we give the graphs of  $\Phi_t$, for $\rho=1$ (resp. $\rho=-1$) and $t=3/2$, $2$ and $5/2$ (resp. $t=1/4$, $1/2$ and $1$) from top to bottom. In both figures, the dashed lines correspond to the sine function, relevant for $X>1$. Moreover, in the negative curvature case, the dotted (lowest) line corresponds to the limit case, obtained when $t$ goes to infinity (with $\rho$ fixed).


\section{Proof of Theorem~\ref{theo-main}}
\label{sec-preuve}

We first assume that the Markov diffusion triple satisfies a $CD(\rho,n)$ condition and prove properties $i$ and $ii$. 

Let then  $f\in\mathcal A$ be a fixed positive function, which we assume to be larger than an $\epsilon>0$, without loss of generality. Let also $t>0$ be fixed, and define $\Lambda(s)=P_s(P_{t-s}f \log P_{t-s}f)$ for $s \in [0,t]$, so that $\ent{P_t}{f} = \Lambda (t) - \Lambda (0)$. Classical properties of the derivative of $\Lambda$, see for instance~\cite[Section~6.7.2]{bgl-book}, ensure that 
$$
\Lambda'(s)=P_s\PAR{\frac{\Gamma(g)}{g}} = P_s \big( g \, \Gamma(\log g) \big)\quad{\rm and}\quad \Lambda''(s)=2P_s\PAR{g \, {\Gamma_2(\log g)}}
$$
for all $s \in [0,t]$, where $g=P_{t-s}f$. Applying the $CD(\rho,n)$ condition, we obtain 
$$
\Lambda''(s)\geq 2\rho\Lambda'(s) + \frac{2}{n}P_s\PAR{g[L(\log g)]^2} . 
$$
Since $\displaystyle L(\log g)=\frac{Lg}{g}-\frac{\Gamma(g)}{g^2}$ by diffusion property, the Cauchy-Schwarz inequality further gives 
$$
\Lambda''(s)\geq 2\rho\Lambda'(s) +  \frac{2}{nP_tf}\PAR{LP_tf-\Lambda'(s)}^2. 
$$
This computation has been performed for example in~\cite{bakry-ledoux} and \cite[Section 6.7.2]{bgl-book}. 

Let now $\displaystyle a=\frac{2}{nP_tf}, X=\frac{4}{n\rho}\frac{LP_tf}{P_tf}, B=\PAR{\frac{n\rho}{2}  P_tf}^2(X-1)$ and $\displaystyle C=LP_t f-\frac{n\rho}{2}P_tf$. 
The above inequality can be written as 
$$
\Lambda''(s)\geq a\PAR{\PAR{\Lambda'(s)-C}^2+B}, 
$$
or equivalently 
$$
H''(s)\leq -a^2B \, H (s)
$$
for all $s\in[0,t]$, where $H(s)=\exp\PAR{-a(\Lambda(s)-Cs)}$.

Then Lemma~\ref{lem-sturm} below has several consequences. First of all, since $H$ is a positive function, then both cases~$i$ and~$ii$ in the lemma are impossible. This means that necessarily $a^2 B< \frac{\pi^2}{t^2}$, or in other words that for any $t>0$
$$
X< 1+ \frac{\pi^2}{\rho^2 t^2}, 
$$
or equivalently:
$$
\forall t>0,\quad\frac{1}{4n\rho}\frac{LP_t f}{P_tf}< 1+ \frac{\pi^2}{\rho^2 t^2}. 
$$
Since we now know that  $a^2B< \frac{\pi^2}{t^2}$, Lemma~\ref{lem-sturm} again ensures that for all $s \in [0,t]$ 
\begin{equation}
\label{eq-tau}
H(s)\geq \tau_{a^2B}(t-s) H(0)+\tau_{a^2B}(s) H(t). 
\end{equation}

The reverse local logarithmic Sobolev inequality~\eqref{eq-logsob-inverse} is now obtained by taking the derivative in \eqref{eq-tau} at $s=0$, in the form
$$
H'(0)\geq-\tau_{a^2B}'(t) H(0)+\tau_{a^2B}'(0) H(t). 
$$
In the first case where $a^2B\in(0,\frac{\pi^2}{t^2})$, or equivalently $X\in(1,1+\frac{\pi^2}{t^2 \rho^2})$, then this inequality is exactly~\eqref{eq-logsob-inverse} with the sine function when written in terms of $\Lambda$, and then of $P_t f$. In the other case where $a^2B < 0$, that is $X < 1$, then we obtain~\eqref{eq-logsob-inverse} with the hyperbolic sine function. The limit case $a^2B=0$, or equivalently $X=0$, is the limit of any of the first two cases.  Together, the obtained estimates are exactly inequality~\eqref{eq-logsob-inverse}. 

The local logarithmic Sobolev inequality~\eqref{eq-logsob} is obtained in the same way, by taking the derivative of~\eqref{eq-tau} at $s=t$.

\medskip

We now prove the converse implication, namely that $ii$ implies $i$, the case $iii$ being handled by the same method. Observing that inequality~\eqref{eq-logsob-inverse} is an equality when $f=1$, the idea is to let $f=1+\varepsilon h$ with $h\in\mathcal A$ and to perform a second order Taylor expansion of~\eqref{eq-logsob-inverse} in the parameter $\varepsilon$ tending to $0$. The zeroth and first order terms in $\varepsilon$ vanish, and recalling that $\ent{P_t}{f}=\frac{\varepsilon^2}{2} \Big[ P_t (h^2) - (P_t h)^2 \Big]+o(\varepsilon^2)$, we obtain
$$
P_t (h^2) - (P_t h)^2  \geq \frac{e^{2\rho t}-1}{\rho}\Gamma(P_th) + \frac{e^{2\rho t} -1 - 2 \rho t}{\rho^2} \frac{(LP_th)^2}{n} 
$$
for all $t \geq0 $. Now, as in~\cite{bakry-ledoux} for instance,  a second order Taylor expansion in $t$ tending to $0$, with the zero and first order terms vanishing, gives the $CD(\rho,n)$ condition back. 

\bigskip

In this proof we have used the following elementary but useful lemma, which is proved  in~\cite[Thm~14.28]{villani-book1} for instance.
\begin{elem}
\label{lem-sturm}
Let $t>0$, $\lambda\in\R$ and  $h$ a $C^2$ non-negative function on $[0,t]$ such that $
h''\leq -\lambda  h$ on $[0,t]$. 
Then:
\begin{enumerate}
\item If $\lambda>\frac{\pi^2}{t^2}$ then for all $s\in[0,t]$, $h(s)=0$.
\item If $\lambda=\frac{\pi^2}{t^2}$ then for all $s\in[0,t]$, $h(s)=c\sin\left(\frac{s}{t}\right)$,  for some $c\geq0$.
\item If $\lambda<\frac{\pi^2}{t^2}$ then  for all $s\in[0,t],\quad h(s)\geq \tau_\lambda(t-s) h(0)+\tau_\lambda(s) h(t)$,
where 
$$
\tau_\lambda(s)=
\left\{
\begin{array}{ll}
\disp\frac{\sin(\sqrt{\lambda}s)}{\sin(\sqrt{\lambda}t)}& {\rm if} \quad \lambda>0\\
\disp{s}/{t} & {\rm if} \quad \lambda=0 \\
\disp\frac{\sinh(\sqrt{-\lambda}s)}{\sin(\sqrt{-\lambda}t)} &{\rm if} \quad \lambda<0
\end{array}
\right.
$$
\end{enumerate}
\end{elem}

\section{Applications to ultracontractive estimates and Harnack inequalities}
\label{sec-applications}
\subsection{Ultracontractive estimates in positive curvature}
\label{sec-ultra}

In this section we assume that $\rho>0$. We will use the new Li-Yau inequality~\eqref{eq-li-yau-gene}  in force to obtain uniform bounds on $L\log P_tf$, for positive $f$, and then ultracontractive estimates on the semigroup. 

\medskip

Under a curvature-dimension condition $CD(\rho,n)$ with $\rho >0$ and $n\geq 1$, it is classical that a reversible Markov semigroup $(P_t)_{t \geq 0}$ is ultracontractive, that is, $P_tf$ is bounded for all $t>0$ and integrable $f$.

 A way of proving an ultracontractive bound is indeed the following: First, a $CD(\rho,n)$ condition with $\rho>0$  implies a (Nash-type) logarithmic entropy-energy inequality:  
 \begin{equation}
\label{eq-ee}
\ent{\mu}{f^2}\leq \frac{n}{2}\log\PAR{1+\frac{4}{\rho n}\int \Gamma(f)d\mu}
\end{equation}
for any $f\in \mathcal A$ such that $\int f^2d\mu=1$. 
Let us observe that, in this ergodic case where $P_t f$ converges to $\int f \, d\mu$ for large time, then~\eqref{eq-ee} can be recovered by letting $t$ go to infinity in our local inequality~\eqref{eq-logsob} (written for $f^2$ instead of~$f$). Then~\eqref{eq-ee} implies the following ultracontractive estimate: there exists a constant $C$ such~that
\begin{equation}
\label{eq-ultra1}
||P_tf||_{\infty}\leq \frac{C}{t^{n/2}}\int |f|d\mu,\quad t\in (0,1]
\end{equation}
for any $f\in\mathcal A$. These results can be found for instance in~\cite[Chapter 6]{bgl-book}. 

\bigskip

This bound is in fact included in the Li-Yau inequality~\eqref{eq-li-yau-gene} which also gives a quantitative estimate for large time.  Observe indeed that 
$$\Gamma(P_tf)/(P_t f)^2\geq0$$ for any positive $f$ and $t$, so~\eqref{eq-li-yau-gene} gives 
$$
\Phi_t\PAR{\frac{4}{n\rho}\frac{L P_tf}{P_t f}} > 0. 
$$
Let now $\rho, t > 0$ be fixed. By definition \eqref{eq-defgh}, it holds $\Phi_t(0)>0$, 
$$
\lim_{x\rightarrow -\infty}\Phi_t(x)=-\infty\quad {\rm  and}\quad \lim_{x\rightarrow 1+ \frac{\pi^2}{\rho^2 t^2}}\Phi_t(x)=-\infty.
$$
Hence the continuous and strictly concave function $\Phi_t$ admits exactly two roots $\xi_1^t<0<\xi_2^t,$ and is positive in-between and negative outside its roots. Of course $\xi_1^t$ and $\xi_2^t$ depend on $t$ and $\rho$ but not on the dimension $n$. In particular
\begin{equation}\label{eq-ultra}
 \xi_1^t < \frac{4}{n\rho} \frac{L P_t f}{P_t f} < \xi_2^t  
 \end{equation}
 for all  positive $f$ and $t$.

 A first simple consequence is the following: if $t \geq 2/ \rho$, then $\Phi_t(1) = 1/ t -\rho/2$ is non-positive, so necessarily $\xi_2^t \leq 1$ (see also the graph of the second function in Figure~\ref{fig-1}). In other words:
\begin{ecor}[\cite{bakry-qian}]
Assume a $CD(\rho, n)$ condition with $\rho>0$. Then for all positive function $f$ in $\mathcal A$ 
\begin{equation}
\label{eq-bq-easy}
\frac{4}{n\rho}\frac{LP_tf }{P_t f} < 1, \qquad t \geq \frac{2}{\rho}\cdot
\end{equation}
\end{ecor}

Let us remark that~\eqref{eq-bq-easy} provides additional information to the bound~\eqref{eq-esti2} in Corollary~\ref{cor-esti2}. Moreover, inequalities~\eqref{eq-esti2} and~\eqref{eq-bq-easy} and are of course not optimal for large $t$ since $LP_tf$ converges to $0$ when $t$ goes to $+\infty$. Next proposition gives an answer to this issue: it makes quantitative the fact that $\xi^t_1, \xi^t_2 \rightarrow 0$ when $t$ goes to $+ \infty$ (and  $\xi^t_1\rightarrow -\infty, \xi^t_2\rightarrow + \infty$ when $t$ goes to $0$), giving by~\eqref{eq-ultra} corresponding explicit upper and lower bounds for $LP_t f/P_tf$.

\begin{elem}
\label{prop-esti}
Assume that $\rho>0$. Then the roots $\xi^t_1$ and $\xi^t_2$ of $\Phi_t$ are such that
$$
-4 e^{-\rho t} - 4 e^{-2 \rho t} \leq \xi^t_1 \leq -4 e^{-\rho t} + 8 \rho t e^{-2 \rho t}, \qquad t \geq \frac{1}{2 \rho}
$$
$$
4 e^{-\rho t} - 4 e^{-2 \rho t} \leq \xi^t_2 \leq 4 e^{-\rho t} + 8 \rho t e^{-2 \rho t}, \qquad t \geq \frac{6}{\rho}
$$
$$
\xi^t_1 =-\frac{2}{\rho t}+\bigO{1} \quad{\rm and }\quad  \xi^t_2=\frac{\pi^2}{\rho^2 t^2}-\frac{4}{\rho t}+\bigO{1}, \qquad t \to 0.
$$
\end{elem}

\medskip
Observe the compatibility of this last bound with Corollary \ref{eq-esti2}.

\medskip

\begin{eproof}
We first consider the large time bounds on the negative root $\xi^t_1$. For $x<0$, we observe that $\xi^t_1 \geq x$ if and only if $\Phi_t (x) \leq 0$, or if and only if $u=\sqrt{1-x}$, which then is larger than $1$, satisfies  
$$
(u-1) e^{\rho t u} \geq u+1. 
$$
For $u = 1+2 e^{-\rho t}$, and by the elementary inequality $e^v \geq 1+v$, this property holds as soon as $t  \geq1/(2\rho)$. Therefore $\xi^t_1\geq 1-u^2$ for these $t$. For $u=1+2e^{-\rho t} - 4 \rho t e^{-2 \rho t}$, the reverse inequality holds also as soon as $t  \geq1/(2\rho)$. Therefore $\xi^t_1 \leq 1-u^2$ for these $t$. Together, this gives the large time lower and upper bounds on~$\xi^t_1.$

\smallskip

Then we recall that the positive root $\xi^t_2$ belongs to $(0,1]$ as soon as $t \geq 2/\rho$, which we assume. Then, now, $\xi^t_2 \leq x$ if and only if $\Phi_t (x) \leq 0$, or if and only if $(1-u) e^{\rho t u} \geq 1+u$. For $u=1-2e^{-\rho t}$ the reverse inequality holds as soon as $t \geq 2/\rho$, so that $\xi^2_t \geq 1-u^2$ for these $t$. For $u=1-2 e^{-\rho t} - 4 \rho t e^{-2 \rho t}$ the inequality holds as soon as $t \geq 6/\rho$, so that $\xi^2_t \leq 1-u^2$ for these $t$. This gives the large time lower and upper bounds on~$\xi^t_2.$

\smallskip

We proceed in the same manner to get the short time estimates. 
\end{eproof}

\medskip

Let us now see how to turn these estimates on $\xi^t_1$, $\xi^t_2$ into ultracontractivity estimates on the semigroup. 
Given $0<t<s$, integrating the pointwise bound~\eqref{eq-ultra} over the interval $[t,s]$ gives  
$$
\exp\PAR{-\frac{n\rho}{4}\int_t^s \xi_2^udu}\leq \frac{P_t f}{P_s f}\leq \exp\PAR{-\frac{n\rho}{4}\int_t^s \xi_1^udu}.  
 $$
 By ergodicity, letting $s$ go to infinity  implies :
 $$
\exp\PAR{-\frac{n\rho}{4}\int_t^\infty \xi_2^u du}\leq \frac{P_t f}{\int fd\mu}\leq \exp\PAR{-\frac{n\rho}{4}\int_t^\infty\xi_1^u du}, \qquad t>0. 
$$
By Proposition~\ref{prop-esti}, it follows that
$$
\frac{P_t f}{\int f \, d\mu} \leq \exp \Big[ n \big(e^{-\rho t} + \frac{e^{-2 \rho t}}{2} \big) \Big], \qquad t \geq \frac{1}{2 \rho}
$$
and
$$
\frac{P_t f}{\int f \, d\mu} \geq \exp \Big[ -n \big(e^{-\rho t} + \frac{1+ 2 \rho t}{2} e^{-2 \rho t} \big) \Big], \qquad t \geq \frac{6}{\rho} \cdot
$$
Proceeding likewise for the short time estimates finally gives :

\begin{ecor}
\label{cor-conv}
Assume a $CD(\rho, n)$ condition with $\rho>0$ with $n\geq 1$. 
Then there exist constants $C$ and $D$, depending on $\rho$ and $n$, such that  for any positive $f$
\begin{equation}
\label{eq-ultra2}
\Big|P_t f-\int fd\mu\Big|\leq C e^{-\rho t} \int fd\mu , \qquad t\geq 1,
\end{equation}
and 
\begin{equation}
\label{eq-ultra22}
\frac{1}{D t^n} \exp\PAR{-\frac{n\pi^2}{\rho t}}\leq \frac{P_t f}{\int fd\mu}\leq \frac{D}{t^{n/2}}, \qquad t\in (0,1].
\end{equation}
\end{ecor}

\smallskip

By Proposition~\ref{prop-esti}, the constant $C$ in~\eqref{eq-ultra2} can be made explicit in $n$ and $\rho$. Inequality~\eqref{eq-ultra2} does not give the asymptotic behavior of $P_t f$ in a satisfactory way. Indeed it is classical that the $CD(\rho, n)$ condition on a symmetric Markov semigroup with $\rho>0$ and $n>1$ implies a (Poincar\'e) spectral gap inequality with constant $(n-1)/(\rho n)$; it follows that $P_t f$ converges to its mean with an exponential speed with rate $\frac{2 \rho n}{n-1}$ :
$$
\int \PAR{P_t f-\int fd\mu}^2d\mu\leq e^{-\frac{2\rho n}{n-1} t}\int \PAR{f-\int fd\mu}^2d\mu. 
$$
The right rate of convergence has been lost in~\eqref{eq-ultra2}. Let us in fact observe that the rate $\rho$ in~\eqref{eq-ultra2} cannot be improved into $\rho n / (n-1)$ by our method since neither $\xi_1^t$ nor $\xi_2^t$ depend on $n$. 

\medskip

On the other hand, since $\vert P_t f \vert \leq P_t \vert f \vert$, the upper bound in~\eqref{eq-ultra22} can be extended to any $f$, recovering the classical ultracontractivity property~\eqref{eq-ultra1}, together with an explicit lower bound on $P_t f$ which is not included in~\eqref{eq-ultra1}. 

\medskip

Moreover the method is powerful since we have not used the symmetric (self-adjoint) property of the equation.  The strategy  can also be applied in the non-symmetric case to derive  new ultracontractive bounds for more general models.

\bigskip

One can also derive similar bounds on the gradients. For instance: 
\begin{ecor} 
Assume a $CD(\rho, n)$ condition with $\rho>0$ with $n\geq 1$. 
Then there exists an explicit constant $C$, independent of $\rho$ and $n$, such that  for any positive $f$
 $$
\Gamma(\log P_t f)=\frac{\Gamma(P_t f)}{(P_tf)^2}\leq  C n \rho \, e^{-2\rho t}, \qquad t\geq \frac{6}{\rho} \cdot 
$$
\end{ecor}

\begin{eproof}
We use the above notation, together with $s = \rho t$, and assume that $s \geq 6$. Then, by Lemma~\ref{prop-esti}, $\vert X \vert \leq 6 e^{-s} \leq 1$, so that in particular $\sqrt{1-X} \geq 1-6 e^{-s}.$ Therefore, by this bound and the elementary $e^s \geq s^3/6$,
$$
e^{-2s \sqrt{1-X}} \leq e^{12s e^{-s}-2s} \leq e^{2-2s}.
$$
In particular it is smaller than $1/3$, so
$$
\coth (s \sqrt{1-X}) = \frac{1+ e^{-2s \sqrt{1-X}}}{1-e^{-2s \sqrt{1-X}}} \leq 1 + 3 e^{-2s \sqrt{1-X}}.
$$
It follows that
\begin{multline*}
\frac{2}{\rho} \, \Phi_t(X) 
\leq 
X - 2 + 2 \sqrt{1-X} (1+3 e^{-2s \sqrt{1-X}})\\ 
\leq 
X- 2 + 2  \big (1- \frac{X}{2} \big) (1 + 3 e^{2-2s})
\leq 
6 e^{2-2s}.
\end{multline*}
This concludes the argument by Corollary~\ref{cor-li-yau-gene}.
\end{eproof}

\subsection{Estimates in non-positive curvature}
\label{sec-negative-curvature}
Let $\rho < 0$ be fixed. For any $t>0$, by definition~\eqref{eq-defgh}, it holds $\Phi_t(1) = 1/t - \rho/2 >0$,
$$
\lim_{x\rightarrow -\infty}\Phi_t(x)=+\infty\quad {\rm  and}\quad \lim_{x\rightarrow 1+ \frac{\pi^2}{\rho^2 t^2}}\Phi_t(x)=-\infty.
$$
Hence the continuous and strictly concave function $\Phi_t$ admits exactly one root $1<\xi^t<1+ \frac{\pi^2}{\rho^2 t^2}$ (see for instance Figure~\ref{fig-2}). It is positive on $(-\infty, \xi^t)$ and in particular
$$
  \frac{4}{n\rho} \frac{L P_t f}{P_t f}\leq \xi^t  
$$
 for all  positive $f$ and $t$, by the Li-Yau inequality~\eqref{eq-li-yau-gene}. This recovers the bound~\eqref{eq-esti2}. In fact~\eqref{eq-esti2} can not be improved by our method since $\xi^t \sim 1+ \frac{\pi^2}{\rho^2 t^2}$ for large time. Indeed, by direct computation as in the previous section,  
$$
1 + \frac{\pi^2}{\rho^2 t^2} \Big( 1 - \frac{2}{\rho t} \Big) \leq \xi^t \leq 1+ \frac{\pi^2}{\rho^2 t^2}, \qquad t \geq \frac{2}{\vert \rho \vert}  \cdot
$$

\subsection{Harnack inequalities}
\label{sec-harnack}
In this section we assume that the space $E$ is a complete, connected and smooth Riemannian manifold $(M,g)$. This example has been described in Section~\ref{sec-gamma2}. We shall let $d$ denote the Riemannian distance on $M$.

The Li-Yau inequality~\eqref{eq-li-yau-gene} for $\rho \neq 0$ can be written as  
 $$
-\frac{|\nabla P_tf|}{P_t f}\geq -\sqrt{\frac{n}{2} \Phi_{t}\PAR{\frac{4}{n\rho}\frac{LP_tf}{P_tf}}}
$$
where $\Phi_t$ is defined in~\eqref{eq-defgh}, or equivalently, 
 $$
-\frac{|\nabla P_tf|}{P_t f}\geq \Psi_{t,\rho}\PAR{\frac{LP_tf}{P_tf}}
$$
where 
\begin{equation}
\label{eq-psi}
 \Psi_{t,\rho}(x)=-\sqrt{\frac{n}{2}\Phi_t\PAR{\frac{4}{n\rho}x}}. 
\end{equation}
  This extends to $\rho =0$ by letting $\Psi_{t,0}(x)=-\sqrt{n/(2t)+x}$.

\medskip

Basic properties of the function $\Psi_{t,\rho}$ are listed  in the following remark. Their proofs are not complicated, and therefore are omitted.  
\begin{erem}
\label{rem-bijection}
Let $t>0$ be fixed. 
\begin{itemize}
\item Let  $\rho>0$.  The function $\Psi_{t,\rho}$ is defined on the interval 
$$
\mathcal I_{t,\rho}=[{n\rho}\xi_1^t/4,{n\rho}\xi_2^t/4],
$$
 where the roots $\xi_1^t < 0 < \xi_2^t$ of $\Phi_t$ have been defined in Section~\ref{sec-ultra} and depend only on $\rho$ and $t$. Its derivative $\Psi_{t,\rho}'$ is an increasing one-to-one function from $(n\rho\xi_1^t/4 , n\rho\xi_2^t/4 )$ onto $\R$.  The Legendre-Fenchel transform of $\Psi_{t,\rho}$ :  
 $$
 \Psi^*_{t,\rho}(x)=\sup_{y\in\mathcal I_{t,\rho}}\{xy-\Psi_{t,\rho}(y)\}
 $$
  is defined and finite for every  $x \in \R$.  
\item Let $\rho<0$. The function $\Phi_t$ admits only one root $\xi^t>0$ (see Section~\ref{sec-negative-curvature}) and then $\Psi_{t,\rho}$ is defined on the interval 
$$
\mathcal I_{t,\rho}=[n\rho \xi^t/4, + \infty).
$$
 Its derivative $ \Psi_{t,\rho}'$ is an increasing one-to-one function from the interval $(n\rho\xi^t/4,\infty)$ onto $(-\infty,0)$. The Legendre transform of $\Psi_{t,\rho}$ is defined and finite on $(-\infty,0)$. 
\item When $\rho=0$,  the function $\Psi_{t,0}(x)=-\sqrt{n/(2t)+x}$ is defined on $\mathcal I_{t,0}=[-n/(2t),+\infty)$. Its Legendre transform is also defined on $(-\infty,0)$. The case $\rho=0$ appears as limit case of the case $\rho<0$, but not of the case $\rho >0$. 
\item For any $\rho\in\R$ and $t>0$, $\Psi_{t,\rho}$  is strictly convex on its interval of definition. 
\item For any $\rho\in\R$ and $t>0$, $\Psi_{t,\rho}^*$ is non-negative. 
\end{itemize}
\end{erem}

\begin{figure}[!h]
\begin{center}
\includegraphics[width=9cm]{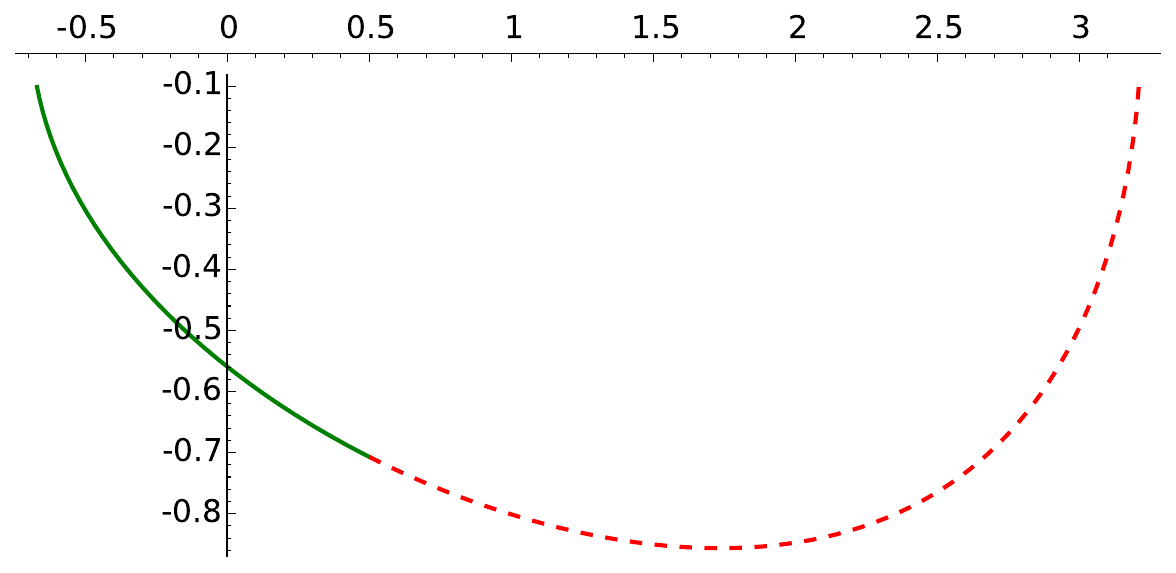}
\caption{Graph of $\Psi_{t,\rho}$ with  a positive  curvature}
\label{fig-3}
\end{center}
\end{figure}

\begin{figure}[!h]
\begin{center}
\includegraphics[width=9cm]{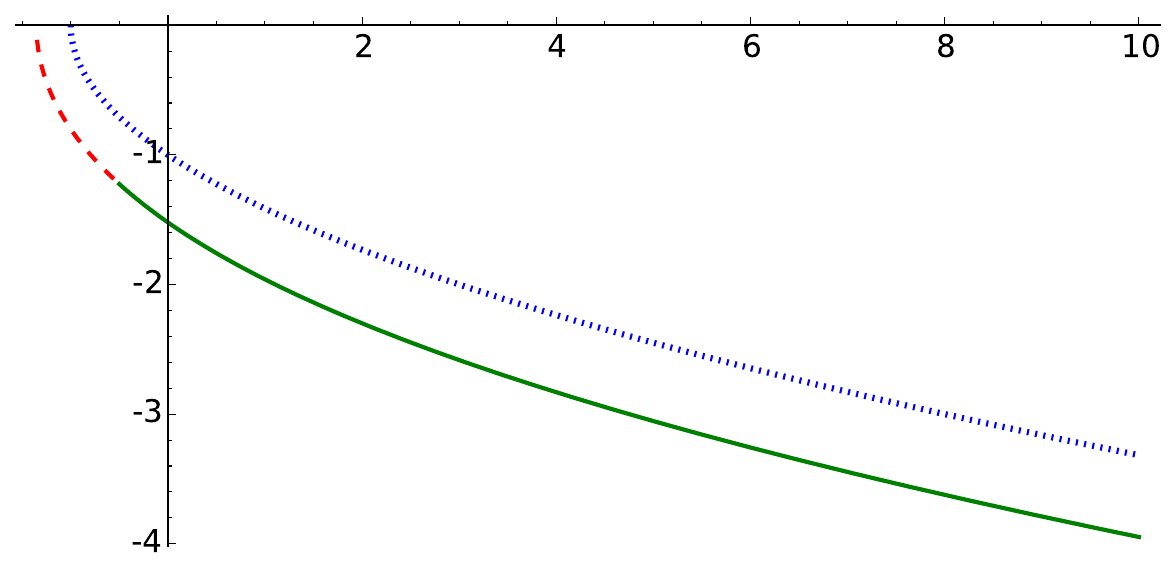}
\caption{Graph of $\Psi_{t,\rho}$ with a negative  curvature}
\label{fig-4}
\end{center}
\end{figure}

In  the figure~\ref{fig-3} (resp. figure~\ref{fig-4}) we have drawn the graph of $x\mapsto\Psi_{t,\rho}(x)$ with  $\rho=1$  (resp. $\rho=-1$),  $t=1$ and $n=2$. The dashed line in the figure~\ref{fig-4} corresponds to the graph of $\Psi_{t,0}$,   again with $t=1$ and $n=2$. 


\bigskip

It is standard since~\cite{li-yau} that a Li-Yau inequality implies a parabolic Harnack inequality in the semigroup. Here we give the {\it equivalence} between these two types of bounds in our framework. 

\begin{ethm}[Harnack inequality]
\label{thm-harnack}
Assume that $L$ satisfies  a $CD(\rho,n)$ condition on the manifold $M$, with $\rho\in\R$ and $n\geq1$.  
\begin{enumerate}
\item Let us assume that  $\rho>0$. 

For  any $s,t>0$,  $x,y\in M$ and any positive function $f$ in $\mathcal A$, we have 
\begin{equation}
\label{eq-harnack}
P_{s}f(x)\leq P_{t}f(y)\exp\PAR{\frac{d(x,y)}{t-s}\int_s^{t}\Psi^*_{u,\rho}\PAR{-\frac{t-s}{d(x,y)}}du}.
\end{equation}

Conversely, if~\eqref{eq-harnack} is satisfied for any positive function $f$ in $\mathcal A$, $x,y\in M$ and $s,t>0$  then the Li-Yau inequality~\eqref{eq-li-yau-gene} holds.  
\item Let us assume that  $\rho\leq0$. 

For  any $0<s<t$,  $x,y\in M$ and any positive function $f$ in $\mathcal A$, then the inequality~\eqref{eq-harnack} holds. 

Conversely, if~\eqref{eq-harnack} is satisfied for any positive function $f$ in $\mathcal A$, $x,y\in M$ and $0<s<t$,  then the Li-Yau inequality~\eqref{eq-li-yau-gene} holds.  
\end{enumerate} 

In both cases \eqref{eq-harnack} also holds for any positive function $f$, in $L^1(\mu)$ for instance.
\end{ethm}

\begin{eproof}
Let us prove the first part of the theorem, when $\rho>0$. 
Let $s,t>0$,  $x,y\in M$ and $f$ positive. Let for $u\in[0,1]$, 
$$
\Lambda(u)=\log P_{b(u)}f(x_u),
$$
where $b(u)=(1-u)s+ut$ ($b$ is not necessarily increasing) and $(x_u)_{u\in[0,1]}$ is a constant speed geodesic between $x$ and $y$. Then we get 
\begin{multline*}
\Lambda'(u)=
(t-s)\frac{LP_{b(u)}f \, (x_u)}{P_{b(u)} f (x_u)}+\frac{\nabla P_{b(u)} f (x_u)}{P_{b(u)}f (x_u)}\cdot \dot{x}_u
\geq\\
(t-s)\frac{LP_{b(u)}f \, (x_u)}{P_{b(u)} f (x_u)} - d \, \frac{\vert \nabla P_{b(u)} f (x_u) \vert}{P_{b(u)}f (x_u)}
\end{multline*}
where $d=d(x,y)$.
The Li-Yau inequality~\eqref{eq-li-yau-gene} ensures  that 
$$
-\frac{|\nabla P_{b(u)} f|}{P_{b(u)} f}\geq \Psi_{b(u),\rho}\left(  \frac{LP_{b(u)} f}{P_{b(u)} f}\right)
$$
at the point $x_u$, where as in~\eqref{eq-psi} $\Psi_{b(u),\rho}(x)=-\sqrt{\frac{n}{2}\Phi_{b(u)}(\frac{4x}{n\rho})}$. But now $\Psi_{b(u),\rho}$ is convex on $\mathcal I_{b(u),\rho}$, so for any $\alpha\in\mathcal I_{b(u),\rho}$
$$
\Lambda'(u)\geq (t-s)\frac{LP_{b(u)} f}{P_{b(u)}}+d \, \Psi'_{b(u),\rho}(\alpha)\left( \frac{LP_{b(u)} f}{P_{b(u)}f}-\alpha\right)+d \, \Psi_{b(u),\rho}(\alpha).
$$
By Remark~\ref{rem-bijection}, for $\rho >0$ there exists $\alpha \in \mathcal I_{b(u),\rho}$ such that 
$$
d \, \Psi'_{b(u),\rho}(\alpha)= -(t-s).
$$
 Hence, for this $\alpha$,
$$
\Lambda'(u)\geq (t-s)\alpha+d \, \Psi_{b(u),\rho}(\alpha).
$$
It gives after integration over $u\in[0,1]$, 
$$
P_{s}f(x)\leq P_{t}f(y)\exp\PAR{d\int_0^1\PAR{-\frac{t-s}{d}\alpha-\Psi_{b(u),\rho}(\alpha)}du}.
$$
Now, by definition of the Legendre-Fenchel transform and of $b(u)$,  
\begin{multline*}
d\int_0^1\PAR{-\frac{t-s}{d}\alpha-\Psi_{b(u),\rho}(\alpha)}ds
\leq
d\int_0^1\Psi^*_{b(u),\rho}\PAR{-\frac{t-s}{d}}du
=\\
\frac{d}{t-s}\int_s^t \Psi^*_{u,\rho}\PAR{-\frac{t-s}{d}}du.
\end{multline*}
This concludes the argument.

\smallskip

Let us now prove the converse part. Let $s>0$ and $y \in M$ be fixed. Let also $x_\varepsilon$ (for $\varepsilon>0$) be the exponential map starting from $y$ and with initial tangent vector $w$. We apply~\eqref{eq-harnack} with $x=x_\varepsilon$ and $t=s+\varepsilon a $, $a\in\{-1,+1\}$. For $f$ positive in $\mathcal A$, 
 inequality~\eqref{eq-harnack} becomes 
\begin{equation}
\label{eq-h2}
P_{s}f(x_\varepsilon)\leq P_{s+\varepsilon a}f(y) \; \exp\PAR{\frac{\varepsilon |w|+\smallO{\varepsilon}}{\varepsilon a}\int_s^{s+\varepsilon a}\Psi^*_{u}\PAR{-\frac{\varepsilon a}{\varepsilon |w|+\smallO{\varepsilon}}}du}.
\end{equation}
A first-order Taylor expansion of~\eqref{eq-h2}  in $\varepsilon>0$ tending to $0$ gives
$$
\frac{\nabla P_sf}{P_s f}\cdot \frac{w}{|w|}\leq \frac{a}{|w|}\frac{LP_sf}{P_sf}+\Psi_{s,\rho}^*\PAR{-\frac{a}{|w|}}
$$
at the point $y$. For $w=r\frac{\nabla P_sf}{|\nabla P_sf|}$ with $r>0$, this can be written as
$$
\frac{|\nabla P_sf|}{P_s f}\leq \frac{a}{r}\frac{LP_sf}{P_sf}+\Psi_{s,\rho}^*\PAR{-\frac{a}{r}}. 
$$
Hence
$$
\frac{|\nabla P_sf|}{P_s f}\leq -\BRA{\frac{LP_sf}{P_sf}z-\Psi_{t,\rho}^*\PAR{z}}  
$$
for any $z\in\R$, since $r>0$ and $a\in\{-1,+1\}$ are arbitrary. Since $\Psi_{s,\rho}^*$ is a convex function on $\mathbb R$, taking the infimum over $z$ finally gives
$$
\frac{|\nabla P_sf|}{P_s f}\leq -\Psi_{s,\rho}\PAR{\frac{LP_sf}{P_sf}}
$$
at the point $y$. This is the Li-Yau inequality~\eqref{eq-li-yau-gene} at the arbitrary time $s>0$ and point $y \in M$.

\medskip

When  $\rho\leq0$, as explained in Remark~\ref{rem-bijection}, $\Psi_{t,\rho}'$ is a one-to-one function from $(n \rho \xi^t /4,\infty)$ onto $(-\infty,0)$. We can use the same method as in the above case $\rho>0$ but one can find such an $\alpha$ only if $0 < s < t$. In other words the argument works only for increasing functions $b$.  The proof of the converse part is also the same in this case by considering only the $0<s<t$. 
\end{eproof}

\begin{erem}
Theorem~\ref{thm-harnack} could be formally stated in the following general form : given a family $\Psi_{t}$ of convex functions, an inequality of the form 
$$
-\frac{|\nabla P_tf|}{P_t f}\geq \Psi_{t}\PAR{\frac{LP_tf}{P_tf}},
$$
for any function $f>0$ on $M$ and $t>0$ is equivalent to a Harnack inequality
 $$
P_{s}f(x)\leq P_{t}f(y)\exp\PAR{\frac{d(x,y)}{t-s}\int_s^{t}\Psi^*_{u}\PAR{-\frac{t-s}{d(x,y)}}du}
$$
for any $s,t>0$ and $x,y\in M$ such that $-(t-s)/d(x,y)$ is in the domain of the Legendre transforms $\Psi^*_{u}$. 
\end{erem}

\begin{erem}
In the limit case $\rho=0$, we have $\Psi_{t,0}(x)=-\sqrt{\frac{n}{2t}+x}$ for $x\geq-\frac{n}{2t}$. Then, $\Psi^*_{t,0}$ is only defined on $(-\infty,0)$, with  $\Psi^*_{t, 0}(y)=-\frac{ny}{2t}-\frac{1}{4y}$ for $y<0$.  Therefore
\begin{multline*}
\frac{d}{t-s}\int_s^{t}\Psi^*_{u,0}\PAR{-\frac{t-s}{d}}du=\frac{d}{t-s}\int_s^t\PAR{\frac{n(t-s)}{2ud}+\frac{1}{4(t-s)}}du=\\
\frac{n}{2}\log\PAR{\frac{t}{s}}+\frac{d^2}{4(t-s)} 
\end{multline*}
for $0<s<t$. Hence, under the $CD(0,n)$ condition, we recover the classical Harnack inequality
$$
P_sf(x)\leq P_{t}f(y)\PAR{\frac{t}{s}}^{n/2}\exp\PAR{\frac{d(x,y)^2}{4(t-s)}}, \qquad 0<s<t, \; x,y\in M.
$$
This Harnack inequality is equivalent to the classical Li-Yau inequality~\eqref{eq-li-yau}.
\end{erem}

Let us now assume that the Markov semigroup admits a density kernel, that is, a function $p_t(x,y)$ such that for any function $f$, $P_tf(x)=\int f(y) \, p_t(x,y) dy$ where $dy$ is the Riemannian measure. This is for instance the case if the semigroup is ultracontractive, so in particular if $\rho >0$.  Then a Harnack inequality classically implies a bound on the kernel. Here we obtain :

\begin{ecor}[Heat kernel bound]
Under the $CD(\rho, n)$ condition, for any $x,y,z\in M$ it holds
\begin{equation}
\label{eq-heat}
p_{s}(z,x)\leq p_t(z,y)\exp\PAR{\frac{d(x,y)}{t-s}\int_s^{t}\Psi^*_{u,\rho}\PAR{-\frac{t-s}{d(x,y)}}du}
\end{equation}
for all $s,t>0$ if $\rho>0$, and all $0<s<t$ if $\rho\leq0$. 
\end{ecor}

\section{Comparison with earlier bounds}
\label{sec-comparaison}

\subsection{Linearisation of the  Li-Yau inequality}
For all given $t>0$ and $\rho \neq 0$, the function $\Phi_t$ is concave (see Section~\ref{sec-result}). Hence the new Li-Yau inequality~\eqref{eq-li-yau-gene} admits a family of linearisations, which is equivalent to it: 
\begin{eprop}[Linearisation of the Li-Yau inequality]\label{prop-lixualpha}
Under the $CD(\rho,n)$ condition with $\rho\in\R^*$ and $n \geq 1$
\begin{enumerate}
\item for any $\alpha\geq0$, 
$$
\frac{\Gamma(P_tf)}{(P_tf)^2}\leq A_1(\alpha)\frac{LP_tf}{P_tf}+\frac{n}{2}B_1(\alpha), 
$$
with 
$$
 A_1(\alpha)=1-\frac{\rho}{2\alpha\sinh^2(\alpha t)}\PAR{\sinh(2\alpha t)-2\alpha t}
 $$
 and 
\begin{multline*}
B_1(\alpha)=\frac{\alpha}{4\sinh^2(\alpha t)}\PAR{\sinh(2\alpha t)+2\alpha t}-\rho+\\
\frac{\rho^2}{4\alpha\sinh^2(\alpha t)}\PAR{\sinh(2\alpha t)-2\alpha t}.
\end{multline*}
\item for any  $\beta\in(0,\pi/t)$, 
$$
\frac{\Gamma(P_tf)}{(P_tf)^2}\leq A_2(\beta)\frac{LP_tf}{P_tf}+\frac{n}{2}B_2(\beta), 
$$
with 
$$
A_2(\beta)=1-\frac{\rho}{2\beta\sin^2(\beta t)}\PAR{2\beta t-\sin(2\beta t)}
$$
and
$$
B_2(\beta)=\frac{\beta}{4\sin^2(\beta t)}\PAR{\sin(2\beta t)+2\beta t}-\rho+\frac{\rho^2}{4\beta\sin^2(\beta t)}\PAR{2\beta t-\sin(2\beta t)}.
$$
\end{enumerate}
\end{eprop}

\begin{eproof}
By Corollary~\ref{cor-esti2} and concavity of the function $\Phi_t$ on $(-\infty, 1+\frac{\pi^2}{\rho^2 t^2} \big)$,
$$
\frac{\Gamma(P_tf)}{(P_tf)^2} <  \frac{n}{2} \Phi_t \Big( \frac{4}{n \rho} \frac{L P_t f}{P_t f} \Big) \leq \frac{2}{\rho} \Phi'_t(x_0) \, \frac{LP_tf}{P_tf}+\frac{n}{2} \big( \Phi_t(x_0) - x_0 \, \Phi'_t(x_0) \big)
$$
for any $x_0 < 1+\pi^2/(\rho^2 t^2)$. By definition~\eqref{eq-defgh} of $\Phi_t$, choosing $x_0 \leq 1$ gives case $i$ with $\alpha = \rho \sqrt{1-x_0}$; observe then that the bound in $i$ is the same for $\alpha$ and $-\alpha$. Choosing $1 < x_0 < 1+\pi^2/(\rho^2 t^2)$ likewise gives case $ii$ with $\beta = \rho \sqrt{x-1}$.
\end{eproof}

Let us observe that the bounds in Proposition~\ref{prop-lixualpha} can also be recovered by extending the method proposed in~\cite{baudoin-garofalo} (see also~\cite{qian}). In their Proposition 2.4, F. Baudoin and N.~Garofalo use a close semigroup argument to prove that
\begin{multline*}
\frac{\Gamma(P_t f)}{(P_t f)^2}\leq \left(1-2\rho\int_0^tV^2(s)ds \right) \frac{LP_t f}{P_t f}+\\
\frac{n}{2} \left( \int_0^tV'(s)^2 ds+\rho^2\int_0^tV^2(s)ds-\rho \right)
\end{multline*}
for any positive function $V$ on $[0,t]$ such that $V(0)=1$ and $V(t)=0$. They deduce diverse bounds for certain choices of $V$, which turn out to be suboptimal: in fact, it is a direct but rather long computation to see that optimising this inequality with respect to $V$, with $\int_0^tV^2(s)ds$ fixed, gives Proposition~\ref{prop-lixualpha}.  
 
 \medskip
 
If  we choose $\alpha=|\rho|$ in $i.$ of Proposition \ref{prop-lixualpha} we obtain the bound~\eqref{eq:lixu} derived by J.~Li and X.~Xu: 
\begin{ecor}[\cite{li-xu}]
\label{cor-li-xu}
Under the $CD(-K,n)$ condition, then   
\begin{equation}\label{eq-lixu}
\frac{\Gamma(P_tf)}{(P_tf)^2}\leq \left(1+ \frac{\sinh(2 Kt)-2Kt}{2\sinh^2(Kt)}\right) \,  \frac{LP_tf}{P_tf}+\frac{n}{2} K\PAR{1+\coth(Kt)}
 \end{equation}
for all positive $t$ and $f$ in $\mathcal A$.
\end{ecor} 
By the same method as here, the authors in~\cite[Thm 1.3]{li-xu} deduced from~\eqref{eq-lixu} the following Harnack inequality :
for every $0<s<t$ and $f>0$,  
\begin{multline}
\label{eq-lixu-harnack}
P_s f(x) 
\leq
P_t f(y) \left(\frac{e^{2Kt}-1-2Kt}{e^{2Ks}-1-2Ks} \right)^{n/4}\\ 
\exp \left(\frac{d(x,y)^2}{4(t-s)} \Big(1+\frac{t \coth (Kt) - s \coth(Ks)}{t-s} \Big) \right).
\end{multline}
Since our Harnack inequality~\eqref{eq-harnack} is equivalent to our Li-Yau inequality~\eqref{eq-li-yau-gene}, which is stronger than~\eqref{eq-lixu}, it follows that~\eqref{eq-harnack} is also stronger than~\eqref{eq-lixu-harnack}. 

\subsection{Comparison with Davies' estimate \eqref{eq:davies} in negative curvature}

Let us prove that the Li-Yau type inequality \eqref{eq-li-yau-gene} given by Corollary \ref{cor-li-yau-gene} improves upon the bound \eqref{eq:davies} established for $\rho <0$, in the notation $K = -\rho >0$.

In our notation $\Phi_t$ and $X = 4 L P_t f / (n \rho P_t f)$ we have to prove that
$$
\Phi_t(X) \leq -\frac{\alpha K}{2}\, X + \frac{\alpha^2}{t}  + \frac{K \alpha^2}{2(\alpha-1)}
$$
for all $\alpha >1, t>0$ and $X$.

For $X<1$, this can be written as
$$
2\, r \, y \, \coth y \leq \frac{1}{\alpha-1} + 2 \, r \, \alpha^2 + (\alpha-1)  r^2 \, y^2 
$$
for all positive $r = (Kt)^{-1}$ and $y = K t \sqrt{1-X}$. But $y \, \coth y \leq  \sqrt{1+y^2}$ as can be seen by taking squares and using $y \leq sinh (y)$ for $y \geq 0$. Hence
$$
2\, r \, y \, \coth y \leq 2 \, r \sqrt{1+y^2} \leq 2 \, r + 2 \, r \, y \leq 2 \, r + \frac{1}{\alpha-1} + (\alpha-1) \, (ry)^2
$$
by the Young inequality. This proves the claim for $X < 1$.

\bigskip

As regards the case where $X>1$,  we first recall the elementary bound~: $y \, \cot y < 1-y^2/4$ for all $y$ in $(0, \pi)$. Letting indeed $u = y/2 \in (0, \pi/2)$, this is due to
$$
y \, \cot y =u \, \frac{1- \tan^2 u}{\tan u} < u \frac{1-u^2}{u} = 1 - \frac{y^2}{4};
$$
here we use that $\tan u >u$ with $x \mapsto \frac{1-x^2}{x} = \frac{1}{x} -x$ decreasing.

As a consequence, in the notation $r = (Kt)^{-1} >0$ and 
$$
y = K t \sqrt{X-1} \in (0, \pi),
$$
\begin{equation}\label{eq:cotan}
\Phi_t(X) = K - \frac{K}{2} X + K \, r \, y \cot y  \leq \frac{K}{2} \big[ 1+ 2 \, r - r^2 \, y^2 - r \, y^2/2 \big],
\end{equation}
so arranging terms it is enough to prove that the second order polynomial
$$
\frac{1}{\alpha-1}  + 2 \, r (\alpha^2-1) + y^2(r^2(1-\alpha) + r/2)
$$
in $y$ is non-negative.

We now observe that the left-hand side of \eqref{eq-li-yau-gene} is non-negative, hence so is $\Phi_t(X)$, and finally the right-hand side of \eqref{eq:cotan}. In other words necessarily $y^2 \leq 2/r$ (that is, $X \leq 1 + 2/(Kt)$ for any $t>0$). Now the above second order polynomial in $y$ (with zero first order coefficient) is non-negative on the interval $[0, \sqrt{2/r}]$ if and only if it is so at $0$ and $\sqrt{2/r}$, which is the case for any $\alpha >1$ by direct computation. This proves the claim for $X > 1$.

\bigskip

The case $X=1$ is covered as a limit case of any of these two cases, or can be directly considered. This concludes the argument, all cases being covered.

\subsection{Comparison with estimate~\eqref{eq:bakryqian} in negative curvature}

Let us now prove that the Li-Yau type inequality \eqref{eq-li-yau-gene} also improves upon the bound \eqref{eq:bakryqian} established for $\rho <0$, in the notation $K = -\rho >0$.

We first observe that, in the notation $r = (Kt)^{-1}, X = 4 LP_t f / (n\rho)$ and $G = 4 \, \Gamma(P_t f) / (nK)$, estimate~\eqref{eq:bakryqian} can be equivalently formulated as
\begin{equation}\label{eq:BQG}
X \leq -G + 2 \, r + 2  \sqrt{G+1 + 2 r}
\end{equation}

\medskip

In the case where $X<1$, then estimate \eqref{eq-li-yau-gene} can be written as
$$
G \leq 2-X + 2 \, \sqrt{1-X} \, \coth (r^{-1} \sqrt{1-X}).
$$
By the elementary inequality $y \, \coth y \leq \sqrt{1+y^2} \leq 1 + y$, and the Young inequality, this implies
\begin{equation}\label{eq-epsilon}
G \leq 2-X + 2 \, r + \frac{1}{\varepsilon} + \varepsilon (1-X)
\end{equation}
for all $\varepsilon >0$, and then
$$
X \leq -G + 2 \, r + \delta \, G + \frac{1}{\delta}
$$
for all $\delta = \varepsilon/(1 + \varepsilon) \in (0,1)$. If $G >1$ then we take $\delta = 1/\sqrt{G} <1$, giving $X \leq -G + 2 \, r + 2 \sqrt{G}$; if $G \leq 1$ then we let $\delta$ tend to $1$, giving $X \leq -G+2r + G +1 \leq -G +2r + 2 \sqrt{G+1}$. In both cases this improves \eqref{eq:BQG}.

\medskip

In the case where $X>1$, then in the same notation the estimate~\eqref{eq-li-yau-gene} can be written as 
$$
G \leq 2-X + 2 \, \sqrt{X-1} \cot (r^{-1} \sqrt{X - 1}).
$$
Since moreover $y \cot y \leq 1$ for all $y \in (0, \pi)$, this implies the bound 
\begin{equation}\label{eq-BQG2}
X \leq -G + 2 \, r + 2,
\end{equation}
which in turn improves \eqref{eq:BQG}.

\medskip

Since again the case $X=1$ is a limit case or can be treated directly, all cases are covered.

%%%%%%%%%%%%
\subsection{Comparison with Hamilton's estimate~\eqref{eq:hamilton} in negative curvature}

Let us finally prove that~\eqref{eq-li-yau-gene} improves upon the bound~\eqref{eq:hamilton} established for $\rho <0$, in the notation $K = -\rho >0$. Letting $s = Kt, X = 4 LP_t f / (n\rho)$ and $G = 4 \, \Gamma(P_t f) / (nK)$, estimate~\eqref{eq:hamilton} can be written as
\begin{equation}\label{eq:hamiltonG}
X \leq - \, e^{-2s} \, G + \frac{2}{s} \, e^{2s}.
\end{equation}
In the case where $X<1$, then choosing $\varepsilon = 2s$ in \eqref{eq-epsilon} ensures that~\eqref{eq-li-yau-gene} implies the bound
$$
X \leq -\frac{1}{1+2s} G +\frac{5 + 4s + 4 s^2}{2s(1+2s)} \cdot
$$
This in turn implies \eqref{eq:hamiltonG} by direct comparison of both coefficients on the right-hand sides. In the case where $X>1$, then the bound \eqref{eq-BQG2}, where $r = 1/s$, likewise implies \eqref{eq:hamiltonG}. As above this concludes the argument.

\nocite{*}
\bibliographystyle{cdraifplain}
\bibliography{li-yau-biblio}
\end{document}